\documentclass[10pt]{amsart}
\usepackage{amsfonts,latexsym,amsthm,amssymb}
\usepackage[all]{xy}

\newcommand{\Q}{{\mathbb Q}}
\newcommand{\Z}{{\mathbb Z}}
\newcommand{\C}{{\mathbb C}}
\newcommand{\A}{{\mathbb A}}
\newcommand{\R}{{\mathbb R}}
\newcommand{\F}{{\mathcal F}}

\newcommand{\E}{{\mathcal E}}
\newcommand{\co}{{\mathcal O}}

\newcommand{\mh}{{\mathfrak m}_h}
\newcommand{\mg}{{\mathfrak m}_g}
\newcommand{\mone}{{\mathfrak m}_1}

\DeclareMathOperator{\Inv}{Inv}
\DeclareMathOperator{\GL}{GL}
\DeclareMathOperator{\Spec}{Spec}

\DeclareMathOperator{\Td}{Td}
\DeclareMathOperator{\ch}{ch}
\DeclareMathOperator{\Hom}{Hom}
\DeclareMathOperator{\Supp}{Supp}

\newtheorem{thm}{Theorem}
\newtheorem{prop}[thm]{Proposition}

\newtheorem{lemma}[thm]{Lemma}

\theoremstyle{remark}
\newtheorem{remark}[thm]{Remark}
\begin{document}
\numberwithin{thm}{section}
\title[Riemann-Roch for quotients]{Riemann-Roch for quotients and Todd
classes of simplicial toric varieties}
\author{Dan Edidin}
\author{William Graham}
\thanks{Both authors were partially supported by N.S.F. grants.}
\maketitle

Department of Mathematics, University of Missouri, Columbia, MO 65211.

{\em E-mail address:} {\tt edidin\@@math.missouri.edu}

\medskip

Department of Mathematics,
University of Georgia,
Athens, GA 30602

{\em E-mail address:} {\tt wag\@@math.uga.edu}

\medskip

{\em Dedicated to Steven L. Kleiman
on the occasion of his 60th birthday.}
\begin{abstract}
In this paper we give an explicit formula for the
Riemann-Roch map for singular
schemes which are quotients of smooth schemes by diagonalizable
groups. As an application we obtain a simple proof
of a formula for the Todd
class of a simplicial toric variety. An equivariant 
version of this formula was previously
obtained for complete simplicial toric varieties by Brion and Vergne \cite{BV}
using different techniques.
\end{abstract}

\renewcommand{\thefootnote}{\fnsymbol{footnote}}
\footnote[0]{
2000 Mathematics Subject Classification: Primary 14C40; Secondary
19L47, 14L32}
\renewcommand{\thefootnote}{\arabic{footnote}}

\section{Introduction}
Suppose $G$ is a linear algebraic group acting properly (hence with
finite stabilizers) on a scheme $X$, with a geometric quotient Y= $X/G$.  
There is \cite{BFM} a Riemann-Roch map
$\tau_Y:G(Y)_{\Q} \rightarrow CH^*(Y)_{\Q}$, where $G(Y)_\Q$ denotes
the rational Grothendieck group of coherent sheaves on $Y$ and $CH^*(Y)_\Q$ 
denotes
the rational Chow groups of $Y$.
Similarly, there is an equivariant Riemann-Roch map
$\tau^G_X: G^G(X)_{\Q} \rightarrow CH^*_G(X)_{\Q}$, where $ G^G(X)$ is
the Grothendieck group of $G$-equivariant coherent sheaves on $X$, and
$CH^*_G(X)$ denotes the equivariant Chow groups of $X$.  

The purpose of this paper is to obtain a formula (Theorem
\ref{thm.eqrr}) for $\tau_Y$ in the case where $X$ is smooth, $G$ is a
diagonalizable group (for example, a torus), and the ground field is
$\C$.  This is of interest because Riemann-Roch maps for non-smooth
schemes are in general difficult to compute.  The formula we give
implies Kawasaki's Riemann-Roch formula for orbifolds which are global
quotients of smooth spaces by diagonalizable groups.  For
quasi-projective quotients of quasi-projective varieties, Theorem
\ref{thm.eqrr} can be deduced from Toen's Riemann-Roch theorem for
Deligne-Mumford stacks \cite[Thm. 4.10]{Toen}.  Our proof, however,
does not require quasi-projectivity hypotheses; moreover, it is more
elementary than Toen's proof and completely avoids the use of stacks.
An application of this result is a simple proof of a formula of Brion
and Vergne \cite{BV} for the Todd class of a simplicial toric variety.
Related results for equivariant Todd classes of toric varieties were
also obtained by Guillemin \cite{Gu} and Ginzburg-Guillemin-Karshon
\cite{GGK} using an equivariant version of Kawaski's Riemann-Roch
theorem.

The ingredients of our proof are the localization theorem in
equivariant $K$-theory and the equivariant Riemann-Roch theorem of
our previous paper \cite{EGRR}.  
The localization theorem relates the equivariant
$K$-theory of $X$ to the equivariant $K$-theory of $X^g$, the fixed
locus of $g \in G$.  It is necessary to consider the fixed loci for
all $g \in G$, and this leads to an expression for $\tau_Y$ in terms of
$\tau^G_X$ which is more subtle than might be naively expected.

In particular, there is a natural map
$\Inv: G^G(X) \rightarrow G(Y)$ taking the class of an equivariant
coherent sheaf $\F$ to the sheaf of invariant sections $\F^G$. 
Also, there is an
isomorphism $\phi: CH^*_G(X)_{\Q} \rightarrow CH^*(Y)_{\Q}$, taking
the class of an invariant cycle $Z$ to the class $\frac{1}{e_{Z}}[Z/G]$,
where (in characteristic 0) 
$e_{Z}$ is the order of the stabilizer group in $G$ of a general
point of $Z$.  These maps fit into a diagram
$$\begin{array}{ccc}
K^G(X) & \stackrel{\tau^G_X}  \rightarrow & CH^*_G(X)_{\Q} \\
\Inv \downarrow & & \phi \downarrow \\
G(Y) & \stackrel{\tau_Y} \rightarrow & CH^*(Y)_{\Q}.
\end{array}$$

The naive formula for $\tau_Y$ is the formula (for $\alpha \in G^G(X)$)
\begin{equation} \label{e.naive}
\tau_Y(\alpha^G) = \phi \circ \tau^G_X(\alpha).
\end{equation}
This would be true if the above diagram commuted.  However, in general
it does not, unless the action of $G$ is free.  The correct formula
for $\tau_Y(\alpha^G)$ is a sum of terms coming from the fixed point loci
$X^g$ as $g$ runs over the (finite) set of elements of $G$ such that
$X^g$ is nonempty.  The right hand side of \eqref{e.naive} is the
contribution corresponding to the fixed point set of the identity in
$G$, which is of course $X$.  If $G$ acts freely on $X$ then the $X^g$
are empty for $g \neq 1$, which explains why the naive formula is true
in that case.

In a subsequent paper we will work out a more general formula 
for actions of reductive groups on algebraic spaces
of arbitrary characteristic. We will also give an equivariant
version of the the formula in Theorem \ref{thm.eqrr}.

\section{Background}

\subsection{Group actions and quotients}
In this paper all schemes are assumed to be of finite type and 
separated over the ground field $\C$. 
All algebraic groups are assumed to be linear.  
Because we are working over $\C$, every {\it reductive} group
(i.e., a group with trivial unipotent radical) is {\it linearly
reductive} (i.e., every representation splits as a direct sum
of irreducibles).  A group $G$ is {\it diagonalizable} if
it is isomorphic to a closed subgroup of $D_n$, the 
group of diagonal
matrices in $\GL(n)$ \cite[Prop.~8.4]{Borel}.
Since $\C$ is algebraically
closed, \cite[Prop.~8.7]{Borel} implies that a group
is diagonalizable if and only if 
it is isomorphic to $(\C^*)^n \times F$
for some finite abelian group $F$.

If $G$ is an algebraic group
acting on a scheme $X$ then we say that $G$ acts {\it properly}
if the map $G \times X \to X \times X$ given by $(g,x) \mapsto (x,gx)$
is proper. This implies that the stabilizer of any point is proper,
hence finite. However, a group action can have finite stabilizers
without being proper.

If $G$ acts properly on a scheme $X$ a geometric quotient $Y = X/G$
always exists in the category of algebraic spaces, but in general $Y$
may not be a scheme.  To avoid the technicalities of working with
algebraic spaces, in this paper when we refer to a quotient $Y = X/G$
we assume that the quotient is a scheme. If $G$ is reductive,
the orbits of the $G$ action on $X$ 
are closed\footnote{This is automatic if $G$
acts with finite stabilizers, since all orbits have dimension equal
$\dim \; G$.}, and 
$X$  can be covered by $G$-invariant affines then a geometric
quotient exists by \cite[Proposition 1.9]{GIT}. If $X$ is
normal and $G$ is a torus then Sumihiro's theorem \cite{Sumihiro}
states that $X$ can be covered by $G$-invariant
affines.
The properness of the action on $X$ implies that the quotient is a
separated scheme \cite[Lemma 0.6]{GIT} and that the quotient
map $X \to X/G$ is an affine morphism \cite[Proposition 0.7]{GIT}.

An important fact about geometric quotients is the following.
\begin{prop} \label{prop.quotients}
Suppose that a reductive group $G$ acts properly on a scheme $X$ and $X \to Y$
is a geometric quotient. Let $p \colon X' \to X$ be a finite $G$-equivariant
map. Then:

(1) $G$ acts properly on $X'$.

(2) There is a geometric quotient $X' \to Y'$ and a finite map
$q:Y' \to Y$ such that the following diagram commutes:
$$\begin{array}{ccc}
X' & \stackrel{p} \to & X\\
\downarrow  & & \downarrow\\
Y' & \stackrel{q} \to & Y.\\
\end{array}$$
\end{prop}
\begin{proof}
(1) Since $G$ acts properly on $X$, the map
$G \times X \to X \times X$ given by $(g,x) \mapsto (x,gx)$ is proper.
The fiber product $(G \times X)_{X \times X} (X' \times X')$
is isomorphic to the  product $G \times (X'\times_X X')$
via the map $(g,x,x'_1, x'_2) \mapsto (g, x'_1, g^{-1}x'_2)$.
Hence the map $G \times (X' \times _X X') \to X' \times X'$
given by $(g,x'_1,x'_2) \mapsto (x'_1, gx'_2)$ is proper.
The map $G \times X' \to X' \times X'$ factors
as
$$G \times X' \stackrel{1 \times \Delta_p} \to G \times (X'\times_X X')
\to X' \times X'.$$
The first map is proper (in fact, it is a closed embedding, 
since $p\colon X' \to X$ is separated).  We have shown that 
the second map is proper,
so the composition is as well. 
Thus $G$ acts properly on $X'$. 
(Note that this argument only requires that $p$ be separated.)

Since $G$ acts properly on $X$ the quotient map
$X \to Y$ is affine. Hence,
$X$ is covered by $G$-invariant
affines $X_1, \ldots X_k$ where $X_i = \Spec A_i$. Since $X' \to X$
is finite, 
$X'_1, \ldots , X' _k$ is a cover of $X'$ by $G$-invariant affines,
where $X'_i = p^{-1}(X_i) = \Spec (B_i)$ and the $B_i$'s are finite
$A_i$-algebras. As noted above, \cite[Proposition 1.9]{GIT} implies
that there is a geometric quotient $X' \to Y'$.
Since geometric quotients are
categorical quotients (\cite[Proposition 0.1]{GIT}) the $G$-invariant
map $X' \to X \to Y$ factors through the quotient $X' \to Y'$; i.e.,
the diagram 
$$\begin{array}{ccc}
X' & \stackrel{p} \to & X\\
\downarrow  & & \downarrow\\
Y' & \stackrel{q} \to & Y\\
\end{array}$$
commutes.

It remains to show that the induced
map $Y' \to Y$ is finite. To do this we will work locally
on $Y$ and assume that $Y = \Spec C$.
Since the quotient map $X \to Y$ is affine it follows
that $X = \Spec A$.  Since $\Spec A \to \Spec A^G$ is a
a categorical quotient \cite[Theorem 1.1]{GIT}, and
categorical quotients are unique up to isomorphism 
(cf. \cite[Proposition 0.1]{GIT}), it follows
that $C \simeq A^G$.

The map $p$ is affine (being finite) so $X' = \Spec 
B$, which again implies that $Y' = \Spec B^G$.
Because $p$ is finite, $B$ is a finite $A$-algebra; we
must prove that $B^G$ is a finite $A^G$-algebra.
Since $G$ is reductive, $A^G$ and $B^G$ are both of finite type over
$\C$ \cite[Theorem 1.1]{GIT}. Thus $B^G$ is of finite type over $A^G$. Now
if $b \in B^G$ then, since $B$ is a finite $A$-algebra, $b$
is integral over $A$. Let $p(t) \in A[t]$ be a monic polynomial
for $b$. Let $q(t) = Ep(t)$ where $E\colon A[t] \to A[t]$
is the Reynolds operator \cite[Definition 1.5]{GIT}. (Here $A[t]$
is made into a $G$-module by taking the $G$-action on $t$ to be trivial.)
Since $1 \in A$ is $G$-invariant, $q(t)$ is a monic polynomial
with coefficients in $A^G$. By functoriality of $E$, 
$Eb$ is a root of $q(t)$. But $b \in B^G$ so $Eb = b$. Therefore,
$B^G$ is integral over $A^G$.
\end{proof}

\subsection{Invariants and equivariant $K$-theory}

If $X$ is a scheme with a $G$-action, denote by
$K^G(X)$ the Grothendieck group of $G$-equivariant locally
free sheaves; this is a ring, with the
product given by the tensor product of equivariant vector
bundles.  Let $G^G(X)$ denote the Grothendieck group of $G$-equivariant
coherent sheaves; this is a module for the ring 
$K^G(X)$.  Identify the
representation ring $R(G)$ with $K^G(pt)$; then pullback by
the projection $X \rightarrow pt$ gives a natural map $R(G) \rightarrow
K^G(X)$.  In this way, $K^G(X)$ becomes an $R(G)$-algebra
and $G^G(X)$ an $R(G)$-module.

If $X$ is smooth and $G$ is reductive then
every $G$-equivariant coherent sheaf has a finite resolution by
$G$-equivariant locally free sheaves.  This implies
that $G^G(X) \simeq K^G(X)$ \cite[Theorem 2.18]{Thomason}.

If $X \stackrel{\pi} \to Y$ is a geometric quotient and $\F$
is a $G$-equivariant coherent sheaf on $X$,
write $\F^G$ for the sheaf of ${\mathcal O}_Y$-modules
$(\pi_*({\mathcal F}))^G$.

\begin{prop} \label{prop.tame}
If a reductive group $G$ acts properly on a
scheme $X$ and $X \to Y$ is a geometric quotient, then
the assignment $\F \mapsto \F^G$ is an exact functor
from the category of $G$-equivariant coherent $\co_X$-modules
to the category of coherent ${\mathcal O}_Y$-modules. 
\end{prop}
\begin{proof}
Since the quotient map is affine
we may assume that $X = \Spec\; A$ and $Y = \Spec\;A^G$.

We first show that the assignment $\F \mapsto \F^G$ is
exact. Since taking invariants is left exact for any
group $G$, it suffices to prove that any surjection $M \twoheadrightarrow N$ 
of $G$-equivariant finitely generated $A$-modules restricts to a surjection
$M^G \twoheadrightarrow N^G$ of invariant submodules. 
Any $G$-module is a union of its finite-dimensional $G$-submodules
(\cite[Lemma p.~25]{GIT} so we may reduce to the case that
$M,N$ are finite-dimensional. In this case, since $G$ is reductive,
the map $M^G \to N^G$ is surjective by Schur's lemma.

To show that the assignment takes coherent sheaves to coherent
sheaves, we must show that if $G$ is reductive and
$N$ is a finitely generated $G$-equivariant
$A$-module, then $N^G$ is finitely generated
as an $A^G$-module. This fact is well known;
see \cite[Theorem 16.8]{Grosshans} for a proof.
\end{proof}

\begin{prop} \label{prop.comminv}
Suppose that $G$ acts properly on $X$ and $X' \stackrel{p} \to X$
is a finite $G$-equivariant map. Let $Y' \stackrel{q} \to Y$
be the induced finite map of quotients (Proposition \ref{prop.quotients}).
If $\F$ is a $G$-equivariant coherent sheaf on $X'$ then 
$q_*(\F^G) = (p_*\F)^G$ as coherent sheaves on $Y$.
\end{prop}
\begin{proof} As in the proof of Proposition \ref{prop.tame}
we can reduce to the case $X = \Spec A$, $X' = \Spec B$,
$Y = \Spec A^G$, $Y'= \Spec B^G$ and $\F = \tilde{M}$
where $M$ is a $G$-equivariant $B$-module. In
this case both $q_*(\F^G)$ and $(p_*\F)^G$ correspond to
the $A^G$-module $M^G$, where the $A^G$-module structure is
obtained via $A^G \to B^G \hookrightarrow B$.
\end{proof}

As a consequence of Propositions \ref{prop.tame} and \ref{prop.comminv},
the assignment $\F \mapsto \F^G$ induces a map of Grothendieck
groups $\Inv \colon  G^G(X) \to G(Y)$ which is covariant for
finite  morphisms. If $\alpha \in G^G(X)$,
we will write $\alpha^G$ for $\Inv(\alpha)$.

\subsection{The localization theorem for diagonalizable 
  groups in equivariant $K$-theory} 
For simplicity, throughout the
rest of the paper all Grothendieck and Chow
groups are taken with complex coefficients.

In this section $G$ will denote a diagonalizable group.
Let $R = R(G) \otimes \C$ denote the representation ring of $G$ tensored with
$\C$.  Then $R$ coincides with the coordinate
ring $\co(G)$ of $G$ \cite[Proposition p.119]{Borel}. 
If $h \in G$, let $\mh \subset R$ denote
the corresponding maximal ideal.

Suppose that $G$ acts on $X$.  
If $h \in G$, let $i_h: X^h \hookrightarrow X$
denote the inclusion of the fixed locus $X^h$ of $h$ into $X$.  
The localization theorem says that the proper
push-forward $i_{h*}: G^G(X^h)_{\mh} \to G^G(X)_{\mh}$ is an isomorphism
\cite[Theorem 2.1]{Thoduke}; here $G^G(X)_{\mh}$ is the
localization of $G^G(X)$ at the maximal ideal $\mh$ of $R$.

If $N$ is an equivariant vector bundle of rank $r$, we write
$\lambda_{-1}(N)$ for the class $\sum_{i=0}^{r} (-1)^{i} [\Lambda^{i}N]$
in $K^G(X)$.
If $N = \oplus_{i = 1}^{r} L_i$ splits as a direct sum of equivariant
line bundles then $\lambda_{-1}(N) = \prod_{i = 1}^{r}(1 - [L_i])$.

If $X$ is smooth, then $i_h$ is a regular embedding
(note that $X^h$ need not be connected or equidimensional).  Let
$N_h$ denote the normal bundle to this embedding, and
$N_h^*$ the conormal bundle.
Then \cite[Lemma 3.3]{Thoduke}
says ${\lambda_{-1}}(N_h^*)$ is a unit in the ring $K^G(X^h)
\simeq G^G(X^h)$  
and the inverse to 
$i_{h*}$ is the map
$$\frac{i^*_h}{\lambda_{-1}(N^*_h)} \colon
G^G(X)_{\mh} \to G^G(X^h)_{\mh}.$$ 

\begin{prop} \label{p.decomp}
Suppose that $G$ acts on $X$ with finite stabilizers.  Then the natural map
$$
G^G(X) \rightarrow \oplus_h G^G(X)_{\mh}
$$
is an isomorphism, where the sum is over the finite set of $h \in G$
such that $X^h$ is nonempty.
\end{prop}
\begin{proof}
If $G$ acts with finite stabilizers, the proof of
\cite[Proposition 5.1]{EGRR} shows that there is an ideal
$J \subset R$ with 
$R/J$ supported at the points $h \in G$
for which $X^h$ is nonempty,
such that $J R^G(X) = 0$. Hence 
$$G^G(X) = G^G(X)/J G^G(X) = G^G(X) \otimes_{R} R/J.$$
Now $R/J$ is Artinian, so the map 
$$R \to \oplus_{\mh \in \Spec R/J} (R/J)/\mh^{k_h} = 
\oplus_{\mh \in \Spec R/J} (R/J)_{\mh}$$ is an isomorphism; here $k_h$
is an integer such that
$\mh^{k_h} R/J = \mh^{k_h + 1} R/J$. Tensoring with $G^G(X)$ yields
the proposition.
\end{proof}

Henceforth, if $G$ acts on $X$ with finite stabilizers, we will 
view $G^G(X)_{\mh}$ as a subset of $G_G(X)$.  If $f: G_G(X) \rightarrow V$
is any map from $G^G(X)$ to any vector space $V$, we use the 
same letter $f$ to denote the induced map $G^G(X)_{\mh} \rightarrow V$. 
In particular, we obtain $\Inv:G^G(X)_{\mh} \rightarrow G(Y)$.

\begin{lemma}  \label{l.finite}
Suppose $p \colon X' \rightarrow  X$ is a
finite $G$-equivariant map, where $G$ acts with finite stabilizers.
Then $p_* (G^G(X')_{\mh}) \subset G^G(X)_{\mh}$.
\end{lemma}
\begin{proof}
Suppose that $\alpha \in G^G(X')_{\mh}$.  By definition, this
means that $\alpha$ is in the kernel of
$G^G(X') \to G^G(X')_{\mg}$ for all $g \neq h$.  
This is equivalent to saying that for all $g \neq h$
there exists $r \in R \smallsetminus \mg$ such that
$r \alpha = 0$.  The equivariant projection formula
implies that $r p_*(\alpha) = p_*(r \alpha) = 0$.
Hence $p_*(\alpha) \in G^G(X)_{\mh}$.
\end{proof}

\subsection{Equivariant Chow groups and quotients} \label{ss.eqchow}
If $X$ is an equidimensional scheme we will use the notation
$CH^i_G(X)$ to denote the ``codimension $i$'' equivariant Chow group
$A^G_{\dim X -i}(X)$ defined in \cite{EIT} (tensored
with $\C$).  
If $G$ acts with finite stabilizers
then $CH^i_G(X)$ is generated by classes of $G$-invariant
subvarieties \cite[Prop.~5.2]{EGRR}. We reserve
the notation $A^*_G(X)$ for the $G$-equivariant operational
Chow ring defined in \cite{EIT}. When $X$ is smooth we
can identify $A^*_G(X)$ with $CH^*_G(X)$ \cite[Proposition 4]{EIT}.

Now suppose that $G$ acts properly on $X$, with a geometric
quotient $X \to Y$.  By \cite[Theorem 3]{EIT},
there is an isomorphism
$$
\phi: CH^*_G(X) \rightarrow CH^*(Y),
$$
defined as follows.  Let $W \subset X$ be a $G$-invariant
subvariety with image $Z \subset Y$. Set $e_W$ to be the order of the
stabilizer of the generic point of $W$.  Then
$\phi([W]_G) = \frac{1}{e_W}[Z]$; this map descends to give
the isomorphism of Chow groups.  The proof of this theorem 
(cf.~\cite[Proposition 11(a)]{EIT})
shows that if $p: X' 
\rightarrow X$ is a finite equivariant morphism and
$q: Y' \rightarrow Y$ is the corresponding morphism of
quotients, then $\phi_X \circ p_* = q_* \circ \phi_{X'}$.

\subsection{The equivariant Riemann-Roch map} \label{ss.err}
By \cite{EGRR}, there is an equivariant Riemann-Roch
map $\tau_X^G\colon G^G(X) \to 
\prod_{i=0}^{\infty}CH^i_G(X)$, with the same functorial properties as
the non-equivariant Riemann-Roch map of \cite{BFM}. If
$G$ acts with finite stabilizers, then $CH^i_G(X) = 0$ for
$i > \dim X$ so the Riemann-Roch map has target $CH^*_G(X)$.

\begin{prop} \label{p.errfinite}
Suppose that $G$ is a diagonalizable
group acting on $X$ with finite stabilizers.  Then
$\tau_X^G(G^G(X)_{\mh}) = 0$
for $h \neq 1$, and
$\tau_X^G: G^G(X)_{\mone} \to CH^*_G(X)$
is an isomorphism.
\end{prop}

\begin{proof}
Let $\widehat{G^G(X)}$ denote the completion of
$G^G(X)$ with respect to the ideal $\mone$.
By \cite[Theorem 3.1(a), Theorem 4.1]{EGRR} 
$\tau^G_X : G^G(X) \to CH^*_G(X)$ factors
through an isomorphism $\widehat{G^G(X)} \to CH^*_G(X)$.
The proof of Proposition \ref{p.decomp} shows that
there is an integer $k_1$ such that $\mone^{k_1}G^G(X) =
\mone^{k_1 +1}G^G(X)$. Thus 
$$\widehat{G^G(X)} = G^G(X)/\mone^{k_1}G^G(X) = G^G(X)_{\mone}$$
and the map $G^G(X) \to \widehat {G^G(X)}$ is the projection
$G^G(X) \to G^G(X)_{\mone}$. The proposition follows.
\end{proof}

There is an equivariant Chern character map 
$\ch^G: K^G(X) \rightarrow A^*_G(X)$.
Suppose that $X$ is smooth.  If $X$ is covered
by $G$-invariant affines (for example, if $G$ acts properly
and a geometric quotient $X \to Y$ exists), 
then, identifying $G^G(X)$ with
$K^G(X)$, 
the proof of \cite[Theorem 3.1d]{EGRR} implies that
$\tau^G_X(\alpha) = \ch^G(\alpha) \Td^G(T_X)$;
here $\alpha \in G^G(X)$ and $T_X$ is the tangent bundle of
$X$.\footnote{This formula
does not follow from the {\it statement} of Theorem 3.1d if $G$
is not connected. However, the proof only
requires that if $X$ is a scheme then there is a system of
pairs $(U,V)$ such that 
the mixed spaces $X \times^G U$ are schemes.  This is true if
$X$ is covered by a finite number of invariant affines.}

\begin{prop}
Suppose that $X$ is smooth (so we identify $G^G(X)$ with
$K^G(X)$) and that $G$ acts on $X$ with finite stabilizers.
Let $\gamma \in G^G(X)$.  Suppose that $\gamma$ is a unit
in $G^G(X)_{\mone}$.
Then $\ch^G(\gamma)$ is a unit in $CH^*_G(X)$.  If 
$\frac{\alpha}{\gamma} \in G^G(X)_{\mone}$, then
$$
\tau^G_X(\frac{\alpha}{\gamma}) = \frac{1}{\ch^G(\gamma)} \tau^G_X(\alpha).
$$
\end{prop}
\begin{proof}
By \cite[Theorems 5.1 and 6.1]{EGRR}, there is an isomorphism
$\widehat{G^G(X)} \to G^G(X)_{\mone}$ where
$\widehat{G^G(X)}$ is the completion along the augmentation ideal $I$
of bundles of virtual rank 0
in $K^G(X) = G^G(X)$. Since $\gamma$ is invertible after localizing
at $\mone$, it is invertible in the completion. Hence $\gamma \notin
I G^G(X)$, so it has non-zero virtual rank.
Hence $\ch^G(\gamma) = {\rm rk}\; \gamma + A$ where $A \in CH^{> 0}_G(X)$.
Thus $ch^G(\gamma)$ has an inverse in the completion 
$\prod_i CH^i_G(X)$. However, since $G$ acts with finite
stabilizers, $CH^i_G(X) = 0$ for $i > \dim X$, so the completion
is just $CH^*_G(X)$.  Then
$$
\tau^G_X(\alpha) = \tau^G_X(\frac{\alpha}{\gamma} \cdot \gamma)
= \tau^G_X(\frac{\alpha}{\gamma}) \cdot \ch^G(\gamma).
$$
Multiplying both sides by the inverse of $\ch^G(\gamma)$ gives
the result.
\end{proof} 

\subsection{Subgroup actions on equivariant $K$-theory for diagonalizable
group actions.}
In this section we again assume that $G$ is diagonalizable
and let $R = R(G) \otimes \C$.

Let $H \subset G$ be a finite subgroup and set $\overline{G} = G/H$,
$\overline{R} = R(\overline{G}) \otimes \C$.  Define an $H$-action
on $R$ by
\begin{equation} \label{e.action1}
(h \cdot r)(g) = r(h^{-1}g).
\end{equation}
Then $h \cdot \mg = {\mathfrak m}_{hg}$, and $R^H = \overline{R}$.

Assume that $H$ acts trivially on $X$.  If $\E$ is a $G$-equivariant
sheaf on $X$, and $U \subseteq X$ is any open subset, then
any $h \in H$ induces a map (also denoted $h$)
$$
h: \E(U) \to \E(U).
$$
We can decompose $\E(U)$ into a sum of weight spaces for $H$:
$$
\E(U) = \oplus_{\chi \in \hat{H}} \E_{\chi}(U),
$$
where $\hat{H}$ is the character group of $H$, and
$\E_{\chi}(U)$ is the $\chi$-weight space of $\E(U)$.
This yields a decomposition of $\E$ into $H$-eigensheaves:
$$
\E = \oplus_{\chi \in \hat{H}} \E_{\chi}.
$$
Each $H$-eigensheaf is $G$-equivariant, so we have 
$[\E] = \sum [ \E_{\chi} ]$ in $G^G(X)$.  Define an
$H$-action on $G^G(X)$ by the rule
\begin{equation} \label{e.action2}
h \cdot \sum [ \E_{\chi} ] = \sum \chi(h^{-1}) [ \E_{\chi} ].
\end{equation}
If $X$ is a point, then $G^G(X) = R$ and the
action \eqref{e.action2} agrees with the action \eqref{e.action1}.
For general $X$ (with $H$ acting trivially), the $H$-actions
on $R$ and on $G^G(X)$ are compatible, in the sense that, for
$r \in R$, $\alpha \in G^G(X)$,
$$
h \cdot (r \alpha) = (h \cdot r)(h \cdot \alpha).
$$

This implies that $h \cdot \mg G^G(X) = {\mathfrak m}_{hg} G^G(X)$,
and, in particular, 
$$h^{-1} \cdot \mh G^G(X) = \mone G^G(X).$$
For simplicity of notation, we will sometimes
write $\alpha(h)$ for $h^{-1} \cdot \alpha$.

We remark that there is a decomposition
$$
G^G(X) \simeq G^{\overline{G}}(X) \otimes_{\overline{R}} R,
$$
and the $H$-action on $G^G(X)$ defined in \eqref{e.action2}
coincides with the $H$-action on 
$G^{\overline{G}}(X) \otimes_{\overline{R}} R$
induced by the $H$-action on $R$.

The following lemma will be used in proving the Riemann-Roch
theorem for quotients.

\begin{lemma} \label{l.inv}
Suppose that $G$ acts on $X$, and $H \subset G$ acts trivially on $X$.
If $\alpha \in G^G(X)$, $h \in H$, then $\alpha^G = (h \cdot \alpha)^G$.
\end{lemma}
\begin{proof}
We may assume that $\alpha = [\E] = \sum [ \E_{\chi} ]$.  Then
(with $\chi_1$ denoting the trivial character of $H$),
$[ \E_{\chi} ]^G = 0$ for $\chi \neq \chi_1$, so $[\E]^G
= [\E_{\chi_1}]^G$ and $(h \cdot [\E])^G =
\sum \chi(h) [ \E_{\chi} ]^G = \chi_1(h)[\E_{\chi_1}]^G$.
The result follows since $\chi_1(h) = 1.$
\end{proof}

\section{The Riemann-Roch theorem for quotients of diagonalizable
group actions}
The main theorem of the paper is the following Riemann-Roch theorem
for quotients. 

\begin{thm} \label{thm.eqrr}
Let $G$ be a diagonalizable group acting properly on a
smooth scheme $X$ and let $X \to Y$ be a geometric quotient.
If $h \in G$, let $i_h: X^h \hookrightarrow X$ denote the
inclusion of the fixed locus of $h$, with normal bundle
$N_h$; let $j_h: Y^h = X^h /G \hookrightarrow Y$ denote the
induced inclusion.  
If $\alpha \in K^G(X)$, then 
\begin{equation} \label{e.theorem1}
\tau_Y(\alpha^G) = \sum_{h \in \Supp \alpha}
j_{h*} \circ \phi_{X^h} \circ \tau^G_{X^h} 
(h^{-1} \cdot \frac{i_h^* \alpha}{\lambda_{-1}(N^*_h)}),
\end{equation}
where $\Supp \alpha$ is the (finite) set of points $h\in G$ such that
$\alpha_{h} \neq 0$.
Equivalently,
\begin{equation} \label{e.theorem2}
\tau_Y(\alpha^G) = \sum_{h \in \Supp \alpha} \phi_X \circ i_{h*}\left( 
\frac{\ch^G(i_h^*\alpha(h))}{\ch^G(\lambda_{-1}(N^*_h)(h))} \Td^G(T_{X^h}) \right).
\end{equation}
\end{thm}
\begin{remark}
When $X$ and $Y$ are quasi-projective, Theorem \ref{thm.eqrr} follows
from the Riemann-Roch theorem for Deligne-Mumford stacks proved in
\cite{Toen}. Our proof avoids the quasiprojectivity hypothesis.
Note that the quotient of a quasi-projective
variety need not be quasi-projective. For example, a simplicial
toric variety whose fan is non-polytopal has no ample line bundles.
In particular, our theorem can be applied to compute the Todd class
of such toric varieties (Theorem \ref{thm.todd}).
\end{remark}

\begin{remark}  In the right hand side of \eqref{e.theorem1}, the element
$\frac{i_h^* \alpha}{\lambda_{-1}(N^*_h)}$ is viewed as an element
of $G^G(X^h)_{\mh} \subset G^G(X^h)$.  Since $h^{-1}$ acts trivially
on $X^h$, the results of the previous subsection imply that
$h^{-1} \cdot \frac{i_h^* \alpha}{\lambda_{-1}(N^*_h)}$
is an element of $G^G(X^h)_{\mone} \subset G^G(X^h)$.
\end{remark}

\begin{remark} \label{rem.hinv} 
If $\alpha = [{\mathcal E}]$ where ${\mathcal E}$ is an equivariant
coherent sheaf then the sums in equations \eqref{e.theorem1} and
\eqref{e.theorem2} are in $CH^*(Y) \otimes \Q$ even though the individual
terms are in $CH^*(Y) \otimes \C$.
\end{remark}

The remainder of the section will be devoted to the proof of this theorem.
Observe that \eqref{e.theorem2} follows from \eqref{e.theorem1}
by Section \ref{ss.err}.

Let $\alpha = \sum \alpha_h$ be the decomposition corresponding
to Proposition \ref{p.decomp}.  We will prove that for all $h \in G$,
\begin{equation} \label{e.mainproof1}
\tau_Y((\alpha_h)^G) = 
j_{h*} \circ \phi_{X^h} \circ \tau^G_{X^h} 
(h^{-1} \cdot \frac{i_h^* \alpha}{\lambda_{-1}(N^*_h)}).
\end{equation}

\medskip

{\em Step 1.}  
Suppose that $G$ acts freely on a (not necessarily smooth) scheme $X'$ 
with quotient scheme $Y'$.  We claim that for $\beta \in G^G(X')$,
\begin{equation} \label{e.mainproof2}
\tau_Y(\beta^G) = \phi_X \circ \tau^G_X(\beta).
\end{equation}
Indeed, this follows from \cite[Theorem 3.1(e)]{EGRR},
since we identify
$G^G(X')$ with $G(Y')$ by the map $\beta \mapsto \beta^G$, and 
$CH^*_G(X')$ with $CH^*(Y')$ by the map $\phi_{X'}$.

Note that the freeness of the action implies that
$\beta_h = 0$ for $h \neq 1$;
so $\beta = \beta_1$.  
\medskip

{\em Step 2.}  We prove \eqref{e.mainproof1} for $h = 1$.  In this
case, we want to show that
\begin{equation} \label{e.mainproof3}
\tau_Y((\alpha_1)^G) = \phi_X \circ \tau^G_X(\alpha).
\end{equation}
By Proposition
\ref{p.errfinite}, $\tau^G_X(\alpha_h) = 0$
for $h \neq 1$, so $\tau^G_X(\alpha) = \tau^G_X(\alpha_1)$.

\begin{lemma} 
Let $p \colon X' \rightarrow X$ be a finite surjective $G$-equivariant
morphism of (not necessarily smooth) $G$-schemes, where $G$ 
acts properly on $X$
(and hence also on $X'$ by Proposition \ref{prop.quotients}). 
Then $p_*: G^G(X')_{\mone} \rightarrow G^G(X)_{\mone}$ is
surjective.
\end{lemma}
\begin{proof} By \cite[Theorem 3.1]{EGRR}, we have a commutative diagram
$$
\begin{array}{ccc}
G^G(X')_{\mone} & \stackrel{\tau^G_{X'}}\rightarrow & CH^*_G(X') \\
\downarrow p_*& & \downarrow p_*\\
G^G(X)_{\mone} & \stackrel{\tau^G_{X}}\rightarrow & CH^*_G(X).
\end{array}
$$
The horizontal arrows are isomorphisms by Proposition
\ref{p.errfinite}.  Since $G$ acts properly,
the rational equivariant Chow groups are generated by $G$-invariant
cycles (\cite[Proposition 13(a)]{EIT}), so the right hand 
$p_*$ is surjective.  Therefore
the left hand $p_*$ is surjective. This proves the lemma.
\end{proof}

By \cite[Theorem 6.1]{Seshadri} there exists a finite surjective morphism 
$p \colon X' \rightarrow X$ of $G$-schemes such that
$G$ acts freely on $X'$ (we do not require that
$X'$ be smooth).  By Proposition \ref{prop.quotients} there is 
a geometric quotient $X' \to Y'$
such that the induced map $q: Y' \rightarrow Y$ is also finite.
The preceding lemma implies that there exists $\beta_1
\in G^G(X')_{\mone}$ such that $p_* \beta_1 = \alpha_1$.
By Lemma \ref{l.finite}, $q_*((\beta_1)^G) = \alpha_1^G$.  
Therefore,
$$
\tau_Y((\alpha_1)^G) = \tau_Y(q_*((\beta_1)^G))
= q_* \tau_{Y'}((\beta_1)^G).
$$
By Step 1, this equals $q_* \circ \phi_{X'} \circ \tau^G_{X'}(\beta_1)$.  
As noted in Section \ref{ss.eqchow},
$q_* \circ \phi_{X'} = \phi_X \circ p_*$.  Hence
$$
q_* \circ \phi_{X'} \circ \tau^G_{X'}(\beta_1)
= \phi_X \circ p_* \circ \tau^G_{X'}(\beta_1)
= \phi_X \circ \tau^G_X \circ p_*(\beta_1)
= \phi_X \circ \tau^G_X(\alpha_1).
$$
This completes Step 2.

\medskip

{\em Step 3.}  We prove \eqref{e.mainproof1} for arbitrary $h$.  
As usual, we view $G^G(X)_{\mh}$ and $G^G(X^h)_{\mh}$
as summands in $G^G(X)$ and $G^G(X^h)$, respectively.  Let
$$
\beta_h = \frac{i_h^* \alpha}{\lambda_{-1}(N^*_h)} \in
G^G(X^h)_{\mh}.
$$
By Proposition \ref{prop.comminv}, $(i_{h*}\beta_h)^G = j_{h*}((\beta_h)^G)$.  Hence
\begin{equation} \label{e.mainproof4}
\tau_Y((\alpha_h)^G) = \tau_Y((i_{h*}\beta_h)^G)
= \tau_Y \circ j_{h*}((\beta_h)^G)) =
j_{h*} \circ \tau_{Y^h}((\beta_h)^G),
\end{equation}
Now, $h^{-1} \cdot \beta_h \in G^G(X^h)_{\mone}$,
and by Lemma \ref{l.inv}, 
$(h^{-1} \cdot \beta_h)^G = (\beta_h)^G$.  
By Step 2,
$$
\tau_{Y^h}((h^{-1} \cdot \beta_h)^G) 
= \phi_{X^h} \circ \tau^G_{X^h}(h^{-1} \cdot \beta_h)
= \phi_{X^h} \circ \tau^G_{X^h}(h^{-1} \cdot \frac{i_h^* \alpha}{\lambda_{-1}(N^*_h)}).
$$
Substituting this into
\eqref{e.mainproof4} yields \eqref{e.mainproof1}. This completes
the proof of Theorem \ref{thm.eqrr}.

\section{Todd classes of simplicial toric varieties}
In this section we apply the Riemann-Roch formula of
Theorem \ref{thm.eqrr} to give a formula for the Todd class,
$\tau({\mathcal O_X})$, of a simplicial toric variety $X$.
When $X$ is complete this is the same as the formula of Brion and
Vergne \cite[Theorem 4.1]{BV}.

\subsection{Global coordinates on toric varieties}
We begin by recalling some basic facts about toric varieties,
following notation as in \cite{Fulton} and
\cite{BV}, as well as the construction
of global coordinates in \cite{Cox}.

Let $T$ be a $d$-dimensional torus.
Denote the the lattice of one-parameter subgroups
by $N= \Hom(\C^*, T)$ and the character group by $M = \Hom(T,\C^*)$.
A toric variety $X$ is a normal variety on which $T$ acts,
with a dense orbit isomorphic to $T$.  
Such a variety is determined by its fan
$\Sigma$ in $N_\R$. The toric variety is said to be {\it simplicial}
if every cone in the fan is  generated by 
linearly independent vectors. In this case
the toric variety has finite quotient singularities. If each
cone is generated by elements of part of a basis for $N$, then the toric
variety is smooth. The toric variety $X$ is complete if and only
if the support of $\Sigma$ is $N_\R$.

Denote by $\Sigma(k)$ the set of $k$-dimensional cones of
$\Sigma$.  If $\sigma$ is a cone, let
$\sigma(k)$ denote its $k$-dimensional faces.
A $k$-dimensional cone $\tau \in \Sigma(k)$ determines a $d-k$-dimensional
$T$-invariant subvariety $V(\tau) \subset X$. 

If $\tau \in \Sigma(1)$ is a 1-dimensional cone we will refer to
the primitive lattice vector generating the cone by $n_\tau$. Then
$A_{d-1}(X)$ is generated by $\{[V(\tau)]\}_{\tau \in \Sigma(1)}$.

Define a map $(\C^*)^{\Sigma(1)} \to T$ by $(t_\tau)_{\tau \in \Sigma(1)}
\to \prod_{\tau \in \Sigma(1)} n_\tau(t_\tau)$.
Let $G$ be the kernel of this map of tori. Let $Z(\Sigma)$ be the
subscheme of $\A^{\Sigma(1)}$ defined by the monomial ideal
$B(\Sigma) =  \langle \hat{x}_\sigma : \sigma \in \Sigma \rangle$, where
$\hat{x}_\sigma = \prod_{\tau \notin \sigma(1)} x_\tau$ and
$\{x_\tau\}_{\tau \in \Sigma(1)}$ are coordinates on $\A^{\Sigma(1)}$.
The group $G$ acts  on $W= \A^{\Sigma(1)} - Z(\Sigma)$ by 
the restriction
of the diagonal $(\C^*)^{\Sigma(1)}$ action.  
Let $a_\tau:(\C^*)^{\Sigma(1)} \to \C^*$ be the projection to the
$\tau$-th factor.  We will also refer to the restriction
of this character to any subgroup of $(\C^*)^{\Sigma(1)}$ as $a_\tau$.

\begin{thm} \cite{Cox}
$X$ is a simplicial toric variety if and only if $X$ is a geometric quotient 
$W/G$. 
\end{thm}

Under the quotient map $W \to X$, the coordinate hyperplane $x_\tau$ maps
to the invariant Weil divisor $V(\tau)$.  More generally,
if $\sigma$ is the $k$-dimensional cone generated by 
$\tau_1, \ldots , \tau_k$, then $V(\sigma)$ is the image 
under the quotient map of the 
the linear subspace $W_\sigma$ defined by the ideal
$(x_{\tau_1}, x_{\tau_2},  \cdots ,  x_{\tau_k})$.
Let $G_\sigma$ be the stabilizer of the linear subspace $W_\sigma$.
By definition, 
\begin{eqnarray*} G_\sigma & = &
\{g \in G \ | \ a_\tau(g) =1\;\;\; \mbox{if}\;\; \tau \notin \sigma(1)\}\\ 
& = & \{(t_\tau)_{\tau \in \sigma(1)}| t_\tau \in \C^*, \; 
\prod_{\tau \in \sigma(1)} n_\tau(t_\tau) = 1\}.
\end{eqnarray*}
This implies that $G_\sigma$ is a finite group whose order is
the order of $N_\sigma/\Sigma_{\tau\in \sigma(1)} \Z n_\tau$,
the multiplicity of the cone $\sigma$.
Let $G_\Sigma$ be the union of the subgroups $G_\sigma$ of $G$.

For any $g \in G$, 
set $$
g(1) = \{ \tau \in \Sigma(1) \ \  | \ \ a_\tau(g) \neq 1 \}.
$$
The fixed locus of $g$ in $\A^{\Sigma(1)}$ is the linear
subspace cut out by the ideal generated by $\{x_\tau\}_{\tau \in g(1)}$
Thus, $W^g$ is the intersection of this linear subspace with
the open set $W \subset \A^{\Sigma(1)}$.
Observe that, by definition, the linear space defined by the ideal
$(x_{\tau_1}, x_{\tau_2}, \ldots , x_{\tau_k})$ is in $Z(\Sigma)$
if and only if the rays $\tau_1, \ldots \tau_k$ 
do not all lie in a single cone of the fan $\Sigma$.
This implies that $W^g$ is nonempty if and only if
$g \in G_\sigma$ for some cone 
$\sigma$.  Hence the support of the 
$R$-module $K^G(W)$ is $G_\Sigma \subset G = \Spec R$.

\subsection{Todd classes of simplicial toric varieties}
The Todd class of the toric variety $X$ is by definition
$\Td(X) := \tau({\mathcal O}_X)$.  Since
${\mathcal O}_X = ({\mathcal O}_W)^G$, we can apply
Theorem \ref{thm.eqrr} to obtain a formula for $\Td(X)$.
When $X$ is complete this formula was previously obtained
by Brion and Vergne \cite[Theorem 4.1]{BV}.

For $\tau \in \Sigma(1)$ let $z_\tau \in CH^1_G(W)$ be the equivariant
fundamental class of the invariant divisor
$\{x_\tau = 0\}$.
The generic stabilizer of this divisor is
$G_{\tau}$, which by the preceding subsection is trivial.
Hence $F_\tau := \phi(z_\tau)$
is the fundamental class of the Weil divisor
$V(\tau)$.

\begin{thm} \label{thm.todd}
\begin{equation} \label{e.todd}
\Td(X) = \tau({\mathcal O}_X) = \sum_{g \in G_\Sigma} 
\prod_{\tau \in \Sigma(1)}
\frac{F_\tau}{1 - a_\tau(g)e^{-F_t}}.
\end{equation}

\end{thm}

\begin{proof}
By Theorem \ref{thm.eqrr} 
\begin{equation} \label{e.toddx} 
\Td(X) = \tau({\mathcal O}_W^G) = \sum_{g \in G_\Sigma} \phi \circ i_{g*}
\left( \frac{\Td^G(T_{W^g})}
{\ch^G(\lambda_{-1}(N^*_{W^g}(g))}\right).
\end{equation}

As noted above, $W^g$ is the intersection of $W$ with the linear subspace
cut out by the ideal generated by $\{x_\tau\}_{\tau \in g(1)}$.
Hence
\begin{equation} \label{e.fundclass}
[W^g] = \prod_{\tau \in g(1)} z_\tau \in CH^*_G(W).
\end{equation}
The equivariant normal bundle to $W^g \subset W$ is the pullback
to $W^g$ of the direct sum
of equivariant line bundles
$$
\oplus_{\tau \in g(1)} {\mathcal O}(x_\tau).
$$
Observe that this direct sum decomposition is
also a decomposition into eigenbundles for the finite subgroup 
$H = \langle g \rangle$,
where the weight on the $\tau$-th factor is the character $a_\tau$. 
Since $c_1({\mathcal O}(x_\tau)) = z_\tau$, we obtain
\begin{equation} \label{e.chw} \ch^G(\lambda_{-1}(N^*_{W^ g}(g)) = i_g^*
(\prod_{\tau \in g(1)}
{1 - a^{-1}_\tau(g)e^{-z_\tau}}).
\end{equation}
Likewise, 
the tangent bundle to $W^g$ splits
as
$$
\oplus_{\tau \notin g(1)}{\mathcal O}(x_\tau).
$$
Hence, 
\begin{equation} \Td^G(T_{W^g}) = i_g^{*}(
\prod_{\tau \notin g(1)} \frac{
z_\tau}{1 - e^{-z_\tau}}).
\end{equation} 
If $\tau \notin g(1)$, then $a_\tau(g) = 1$, so
\begin{equation} \label{e.toddw}
\Td^G(T_{W^ g}) = i_g^{*}(\prod_{\tau \notin g(1)} 
\frac{z_\tau}{1 - a_\tau(g^{-1})e^{-z_\tau}}).
\end{equation} 

Substituting \eqref{e.fundclass}--\eqref{e.toddw} into equation \eqref{e.toddx}
shows that the contribution of $W^g$ to the Todd class is
\begin{equation} 
\phi(\prod_{\tau \in \Sigma(1)} \frac{z_\tau}{1 - a_\tau(g^{-1})e^{-z_\tau}})
\end{equation}
Summing over all $g \in G_\Sigma$ and noting that $\phi(z_\tau) = F_\tau$
yields the Brion-Vergne formula. (Note that we are using the fact
that $g \in G_{\Sigma}$ if and only if $g^{-1} \in G_{\Sigma}$.)
\end{proof}

\begin{remark}
When $X$ is complete, then the image of formula \eqref{e.todd} in
cohomology is the same as the restriction from $T$-equivariant 
cohomology of the formula in
\cite[Theorem 4.1]{BV} to ordinary cohomology.
The methods
of this paper can also be used to yield a formula for
the $T$-equivariant Todd class in $CH^*_T(Y)$,but we do not pursue this here.
\end{remark}

\medskip

\noindent{\bf Acknowledgments.} The authors are grateful
to Michel Brion for helpful comments on an earlier version of the paper.

\end{document}